\documentclass[30 pt]{amsart}
\usepackage[latin1]{inputenc}
\usepackage{amssymb}
\usepackage{latexsym}
\usepackage{amscd}
\usepackage{longtable}
\usepackage{amsmath}
\usepackage{array}
\usepackage[all]{xy}
\usepackage{a4}
\newtheorem{teor}{Theorem}[section]

\newtheorem{lema}[teor]{Lemma}

\newcommand{\Q}{\mathbb{Q}}

\newcommand{\Z}{\mathbb{Z}}

\newcommand{\C}{\mathbb{C}}
\newcommand{\R}{\mathbb{R}}

\newcommand{\SL}{\operatorname{SL}}

\newcommand{\Aut}{\operatorname{Aut}}

\newcommand{\Jac}{\operatorname{Jac}}

\newcommand{\ord}{\operatorname{ord}}

\newcommand{\cG}{{\mathcal G}}

\begin{document}

\title[Equations of Bielliptic Modular Curves ]{Equations of Bielliptic Modular Curves }
\author{Josep Gonz\'{a}lez}
\address{ Josep Gonz\'{a}lez,\newline
\indent Universitat Polit\`{e}cnica de Catalunya,\newline \indent
Departament de Matem\`{a}tica Aplicada IV (EPSEVG),\newline \indent
Av.\ Victor Balaguer s/n, \newline \indent 08800 Vilanova i la
Geltr\'{u}, Spain.} \email{josepg@ma4.upc.edu}

\thanks{The  author is partially supported  by DGICYT Grant MTM2009-13060-C02-02.\\
\phantom{ccc}\,2010 {\it Mathematics Subject Classification}: 11F03, 14H45.}

\maketitle

\begin{abstract}
We give a procedure to determine equations for the modular curves $X_0(N)$ which are bielliptic and equations for the $30$ values of $N$ such that $X_0(N)$ is bielliptic and nonhyperelliptic are presented.
\end{abstract}

\section{Introduction}

A  curve $X$ (smooth and projective) of genus $g>1$ defined over a number field $K$ is said to be hyperelliptic, resp. bielliptic, over $K$, if there is an involution $w$ defined over $K$ such that the quotient curve $Y=X/\langle w\rangle$ has genus zero, resp. genus one,  and  $Y(K)\neq \emptyset$. The last condition amounts to saying that $Y$ admits a hyperelliptic model over $K$ or the genus one quotient is an elliptic curve over $K$.

According to Abramovich and Harris (cf. \cite{AH}), we know that for a curve $X/K$ of genus  $>1$ the set of the points $P\in X(\overline{K})$ lying in a quadratic extension of $K$ contains infinitely many points if and only if $X$ is hyperelliptic over $K$ or $X$ is bielliptic over $K$ and the corresponding elliptic quotient curve has rank $\geq 1$.

When we restrict to the modular curves $X_0(N)/\Q$,
the  cusp $\infty\in X_0(N)(\Q)$ and, thus,    $X_0(N)$  is  hyperelliptic or bielliptic over $\Q$ if there exists an involution $w$ defined over $\Q$ whose quotient curve has genus $\leq 1$.
  In  \cite{Ogg}, Ogg  determined the $19$ values of $N$ for which   $X_0(N)$ is hyperelliptic over $\Q$. Later in  \cite{Bars}, Bars determined the $41$  values of $N$ for which $X_0(N)$ is bielliptic over $\Q$. Next, we display the unique   $30$ values of $N$ of all them  for which $X_0(N)$ is non-hyperelliptic:
    \begin{equation}\label{N}
  \begin{array}{c}
  34, 38, 42, 43, 44, 45, 51, 53, 54, 55, 56, 60, \phantom{1}61, \phantom{1}62, \phantom{1}63,\\ [3 pt] 64, 65, 69, 72, 75, 79, 81, 83, 89, 92, 94, 95, 102, 119, 131.
  \end{array}
\end{equation}
  In \cite{Go91}, equations for the $19$ hyperelliptic modular curves are presented. The goal of
  this article is to determine equations for these $30$ bielliptic modular curves. In this way, for each modular curve $X_0(N)$ of genus $>1$ we could determine almost  all  points lying in a quadratic field.

\section{General facts on the modular curve $X_0(N)$}

From now on, $X_0(N)$ denotes the algebraic curve over $\Q$ attached
to the modular group $\Gamma_0(N)$ and such that $\Q(X_0(N))$ is the subfield of $\C(X_0(N))$
which consists of the functions with rational $q$-expansion at the cusp $\infty$, where $q=e^{2\,\pi\,iz}$. Next, we summarize  some well-known facts which shall be used in the next section.

\subsection{The group $\Aut (X_0(N))$}
The group $\SL_2(\R)/\{ \pm 1\}$ is the group of the automorphisms of the complex upper half-plane. Let us denote by $\Gamma_0^*(N)$ the normalizer of $\Gamma_0(N)$ in $\SL_2(\R)/\{ \pm 1\}$. The group $B(N)=\Gamma_0^*(N)/\Gamma_0(N)$ provides a subgroup of
$\Aut (X_0(N))$ described by Lehner and Newman in \cite{LN}  and, later,  revised by Atkin and Lehner in \cite{AL}.

 Let $e_2$ and $e_3$ be the greatest exponents such that $2^{2e2}\cdot 3^{e_3}$ divides $N$. Set $\nu_2=2^{\operatorname{min\,}(3,\lfloor e_2/2\rfloor)}$ and $\nu_3=3^{\operatorname{min\,}(1,\lfloor e_3/2\rfloor)}$. For every positive divisor $d$ of $\nu_2$ or $\nu_3$,  the matrix
$ \left( \begin{array}{cc} 1 & 1/d\\ 0 &1\end{array}\right)\in \Gamma_0^*(N)$ and, thus,  provides an automorphism of $X_0(N)$ which will denoted by $S_d$.

For any positive divisor $d$ of $N$ coprime to $N/d$, the matrix $\frac{1}{\sqrt d}\left(\begin{array}{cc} A\cdot d & B \\ N \cdot C & D\cdot d\end{array}\right)$ with determinant $1$ and $A,B,C,D\in\Z$ lies in $\Gamma_0(N)^*$ and provides an involution
 $w_d$ on $X_0(N)$  independent on $A,B,C$ and $D$, called the  Atkin-Lehner involution attached to  $d $. We denote by $W(N)$ the set of the Atkin-Lehner involutions, which is a commutative group since $w_{d_1}\cdot w_{d_2}=w_{d_1\cdot d_2/\gcd (d_1,d_2)^2}$.  The group $B(N)$ is generated by the group $W(N)$,  $S_{\nu_2}$ and $S_{\nu_3}$. In \cite{KM}, Kenku and Momose proved that when the genus of $X_0(N)$ is $>1$  and  $N\neq 37,63$ one has  that $\Aut (X_0(N))=B(N)$.
\begin{lema}
The field of definition of any Atkin-Lehner involution and $S_2$ is $\Q$ and  for $d> 2$ the field of definition of $S_d$ is  $\Q( \zeta_d)$,   where $\zeta_d$ is a primitive $d$-th root of unity.
\end{lema}

\noindent {\bf Proof.} For an Atkin-Lehner involution $w_d$, one has that $w_d^*(\Q(X_0(N))=\Q(X_0(N))$ and, thus, $w_d$ is defined over $\Q$. Indeed, the  function field  $\Q(X_0(N))$ is generated by the functions $j(z)$, $j(N\,z)$ and
 $w_d$ sends these functions to the functions $j(d\,z)$, $j(N/d\,z)$, which lie in $\Q(X_0(N))$. It is immediate to check that for $S_d$ the number field $\Q( \zeta_d)$ is the smallest number field  $K$ such that $K\otimes \Q(X_0(N))$ contains
 $S_d^*(\Q(X_0(N)))$.
\hfill $\Box$

\subsection{Cusp forms of weight two} We recall that we can identify  the $\C$-vector space  of weight 2 cusp forms on $\Gamma_0(N)$, i.e. $S_2(\Gamma_0(N))$,  with $\Omega^1_{X_0(N)/\C}$ via the map $f(q)\mapsto f(q)\,\frac{dq}{q}$. Moreover, via this map $\Omega^1_{X_0(N)/\Q}$ is in bijective corres\-pondence  to the set of  weight two cusp forms  with rational $q$-expansion.

Let $M$ be a positive divisor of $N$. For any positive divisor $d$ of $N/M$, the map on the complex upper half-plane given by  $z\mapsto d z$ provides a nonconstant morphism $B_d\colon X_0(N) \rightarrow X_0(M)$ which  acts on the cusp forms of weight two
by sending $f(q)\in S_2(\Gamma_0(M))$ to $f(q^d)\in S_2(\Gamma_0(N))$. The vector space $S_2(\Gamma_0(N))^{\operatorname{old}}$ is defined as the sum of the images of such maps for all $M|N$ and $d|N/M$. The vector space  $S_2(\Gamma_0(N))$ has a hermitian inner product called the Petersson inner product and
the vector space $S_2(\Gamma_0(N))^{\operatorname{new}}$ is defined as the orthogonal complement to $S_2(\Gamma_0(N))^{\operatorname{old}}$. We denote by $\operatorname{New}_N$ the set of normalized cusp forms in $S_2(\Gamma_0(N))^{\operatorname{new}}$ which are eigenvectors of all Hecke operators and Atkin-Lehner involutions. By a normalized cusp form we mean a cusp form whose first non-zero Fourier coefficient is equal to $1$. It is well-known that $\operatorname{New}_N$ is a basis of $S_2(\Gamma_0(N))^{\operatorname{new}}$.

In Table 5 of \cite{Taules}, it can be found the dimensions of the vector spaces $S_2(\Gamma_0(N))^G$ and $(S_2(\Gamma_0(N))^{\operatorname{new}})^G$ for any subgroup $G$ of the group of the Atkin-Lehner involutions $W(N)$ for $N\leq 300$. The following result will be useful in order to determine weight two cusp forms invariants under an Atkin-Lehner involution.

\begin{lema}\label{eigenvector}
Let $M$ and $N$ be positive integers such that $ M|N$.  Let  $M_1 $ be a positive divisor of $M$ such that  and $\gcd(M,M/M_1)=1$ and let $\ell$ be  a positive divisor of $N/M$ such that $\gcd(M_1\,\ell,N/(M_1\,\ell)=1$. If $f\in S_2(\Gamma_0(M))$ is a normalized  eigenvector of the Atkin-Lehner involution $w_{M_1}$ with eigenvalue $\varepsilon (f)$ and $\varepsilon\in \{-1,1\}$, then $f(q)+\varepsilon \, f(q^{\ell})\in S_2(\Gamma_0(N))$ is a normalized  eigenvector of the Atkin-Lehner involution $w_{M_1\,\ell}$ with eigenvalue $\varepsilon(f)\cdot \varepsilon$.
\end{lema}
\noindent
{\bf Proof.} An automorphism $u$ on $X_0(N)$ whose action on the  upper half-plane is given by a matrix $\gamma =\left( \begin{array}{cc}A& B\\C&D\end{array}\right)\in \SL_2(\R)$, sends a weight two cusp form $h$ to $h(\gamma(z))\frac{1}{(Cz+D)^2}$. The statement follows from the fact that $w_{M_1\,\ell}$ sends $f$ to $\varepsilon(f)f(q^{\ell})$. \hfill $\Box$
 \subsection{Modular parametrizations of elliptic curves}

Since Shimura-Taniyama-Weil's conjecture was proved, we know that for an elliptic curve $E/\Q$ there exist a positive integer $N$ and a nonconstant morphism $\pi\colon X_0(N)\rightarrow E$ defined over $\Q$. Such a  morphism $\pi$ will be called a modular parametrization of $E$ and level $N$. The following conditions are equivalent:
\begin{itemize}
\item[(i)] The conductor of $E$ is $N$.

\item [(ii)] There exists a modular parametrization $\pi$ of $E$ and level $N$ such that $\pi^*(\Omega^1_E)\subset S_2(\Gamma_0(N))^{\operatorname{new}}\frac{dq}{q}$\,.
\item[(iii)] There exist  $f\in\operatorname{New}_N$ and a modular parametrization $\pi$ of $E$ and level $N$  such that $\pi^*(\Omega^1_{E/\Q})=\Q \cdot f(q)\frac{dq}{q}$.
\end{itemize}
Such a parametrization will be called {\bf new} of level $N$ and the  newform $f$ in part (iii) is unique  and determines  the $\Q$-isogeny class of $E$.

For a non-new modular parametrization $\pi$  of $E$ and level $N$, the conductor $M$ of $E$ divides $N$ and $\pi^*(\Omega^1_{E/\Q})$ is an one dimensional  $\Q$-vector subspace of
$$H_{E,N}=\bigoplus_{ d|N/M} \Q f(q^d)\frac{dq}{q}\subset S_2(\Gamma_0(N))^{\operatorname{old}}\,,$$
where $f$ is the normalized newform of level $M$ attached to $E$.
In fact, for any nonzero cusp form $h\in H_{E,N}$ there exists a modular parametrization $\pi$ of $E$ and level $N$  such that $\pi^*(\Omega^1_{E/\Q})=\Q\cdot h(q)\,\frac{dq}{q}$.

A modular parametrization $\pi$ of $E$ and  level $N$ is called {\bf optimal} if the morphism induced on their jacobians  $\pi_*\colon \Jac(X_0(N))\rightarrow E$ has connected kernel. If $\pi^*(\Omega^1_{E/\Q})=\Q\cdot h(q)\,\frac{dq}{q}$ for some $h\in S_2(\Gamma_0(N))$, then the condition to be $\pi$  optimal amounts to saying that the elliptic curve attached to the lattice
\begin{equation}\label{xarxa}
\Lambda=\left\{\int_{\gamma}h(q)\frac{dq}{q}: \gamma\in H_0(X_0(N),\Z)\right \}
\end{equation} is $\Q$-isomorphic to $E$, i.e. $c_4(\Lambda)=\alpha^4\, c_4(E)$ and $c_6(\Lambda)=\alpha^6\, c_6(E)$  for some nonzero $\alpha\in \Q$. For another modular parametrization $\pi_1$ of an elliptic curve  $E_1/\Q$ and level $N$ such that $\pi^*(\Omega^1_{E/\Q})=\pi_1^*(\Omega^1_{E_1/\Q})$, if $\pi$ is optimal then there is an isogeny $\mu\colon E\rightarrow E_1$ defined over $\Q$ such that $\pi_1=\mu\circ \pi$ and, in particular, $\deg \pi|\deg \pi_1$.

We will denote the $\Q$-isomorphism class of an  elliptic curve $E/\Q$  by giving  its conductor  $N$ and  Cremona's label, i.e. a letter $X$ and a positive integer. For instance, the elliptic curve
$15A8$ stands for  the elliptic curve of conductor $15$ with Cremona's label $A8$. The conductor $N$ and the letter $X$, for instance $15A$, denotes the $\Q$-isogeny class of $E$ and $f_{NX}$ will denote the attached newform to $E$. The optimal quotient in the $\Q$-isogeny class of $E$, called the strong Weil curve,  is always labeled with the number $1$.

We point out  that Manin's conjecture has been checked for all strong Weil elliptic curves, i.e. optimal new modular parametrizations, in Cremona's tables. That is, if  $y^2+a_1 x\,y+ a_3\,y=x^3+a_2\,x^2+a_4\,x+a_6$ is a minimal model over $\Z$ of a strong Weil elliptic curve $E$ of conductor $N$, then there exists a  new nodular parametrization $\pi$ of level $N$ for $E$ such that
$$\pi^*\left (\frac{dx}{2 \,y+a_1\,x+a_3} \right)=\pm f(q)\frac{dq}{q}\,,
$$
where $f\in\operatorname{New}_N$ is the corresponding normalized newform attached to the $\Q$-isogeny class of $E$. In other words, if $\pi(\infty)$ is taken to be the infinity point of $E$, the $q$-expansions of the modular functions $x,y$ are of the form:
 $$
 x=\frac{1}{q^2}+\sum_{n\geq -1} a_n q^n\,,\quad  y=\mp \frac{1}{q^3}+ \sum_{n\geq -2} b_n q^n\quad \text{ and }\quad a_n,b_n\in\Q\,.
 $$
 Equivalently,  when we replace $h$ with $f$  in (\ref{xarxa}),  the lattice $\Lambda$ obtained  and the the minimal model of $E$ have the same invariants $c_4$ and $c_6$.

\section{Procedure to determine equations}
From now on, $N$ is a value in the list (\ref{N}). Let $w$ be a bielliptic involution defined over $\Q$,  let $\pi\colon X_0(N)\rightarrow X_0(N)/\langle w\rangle$ be the natural projection and let us denote by $E$ the elliptic curve  $(X_0(N)/\langle w \rangle,\pi (\infty))$. Since $\deg \pi=2$, the parametrization $\pi$ is optimal. Now we split  in three steps the procedure to find an equation for $X_0(N)$:
\vskip 0.3 cm
\noindent{\it Step 1: Determination of $E$, $w$ and the normalized cusp form $h$ such that $\pi^*(\Omega^1_{E})=\langle  h(q)dq/q\rangle$.}
Since $X_0(N)$ can have several bielliptic involutions defined over $\Q$, first we will determine  for which values of $N$ we can take $w$ such that $\pi$ is new.
Clearly, $X_0(N)$ is bielliptic with a new modular parametrization for an elliptic quotient $E$ if and only if there exists an  elliptic curve of conductor $N$ with modular degree equal to $2$.

By checking in Table 22 of \cite{STNB} among the elliptic curves $E$ with conductor $N$ as in (\ref{N}), we obtain that this fact occurs exactly for $22$ values of $N$. For each of  these  values of $N$, we fix this elliptic curve  $E$ as the co\-rresponding bielliptic quotient and we determine the involution $w\in \Aut_{\Q}(X_0(N))$ such that $(\operatorname{New}_N)^{\langle w\rangle}=h$, where $h$ is the newform attached to $E$ (see Table 1).

For the remaining $8$ values of $N$, i.e. $N\in\{42,60,63,72,75,81,95,119\}$, we fix a bielliptic  involution $w$ defined over $\Q$ among the given in Theorem 3.15 of \cite{Bars}. More precisely, for $N\neq 72$ we choose  $w$ to be an Atkin-Lehner involution and for $N=72$ we take
$w=S_2$  (see Table 2).  In order to   find  a normalized cusp form $h$ such that $S_2(\Gamma_0(N))^{\langle w\rangle}=\langle h\rangle$, we proceed as follows. For $N\neq 72$, by applying  Lemma \ref{eigenvector}, we can easily determine a newform $f\in\operatorname{New}_M$ with  $M|N$ and and integer $q$-expansion such that $S_2(\Gamma_0(N))^{\langle w\rangle}=(\bigoplus_{ d|N/M} \Q f(q^d)\frac{dq}{q})^{\langle w\rangle}=\Q\,h$. The newform $f$ only determines the $\Q$-isogeny class of $E$. In order to determine its $\Q$-isomorphism class, we compute the corresponding lattice $\Lambda$ attached to $h$. This fact allow us to identify $E$ in Cremona's tables  (cf. Table 2). In all cases, $\Lambda$  turns out to be  the lattice corresponding to  a minimal model of the elliptic curve  $E$. For $N=72$, the normalized cusp form $h$ is $f(q^2)$, where $f$ is  the normalized newform  of level $36$ attached to the isogeny class $36A$ and  $E$ is the elliptic curve $36A1$.

\vskip 0.3 cm
\noindent{\it Step 2: Determination of the $q$-expansions of the functions $x,y\in\Q(E)$.} Let $F_N(x,y)=y^2+a_1\,x\,y+a_3\,y-(x^3+a_2\,x^2+a_4\,x+a_6) \in \Z[x,y]$ be the polynomial such that $F(x,y)=0$ is the minimal model given in Cremona's Tables for $E$. For each $N\neq 72$, the lattice obtained from the normalized cusp form $h$ in $\pi^* (\Omega^1_{E/\Q})\subset S_2(\Gamma_0(N))$  corresponds to the  minimal model $F_N(x,y)=0$. Therefore, we can take $y$ such that  the coefficients of the Fourier expansion of $x$ and $y$ are the form
$$
x=\frac{1}{q^2}+\sum_{n\geq -1} a_n q^n\,,\quad y=\frac{1}{q^3}+ \sum_{n\geq -2} b_n q^n\,.
$$
The Fourier coefficients of these modular functions can be determined recursively by means of the equalities
$$
y=-\left(\frac{q d h/dq}{h}+a_1\,x+a_3\right)\left / 2\right.\quad \text{ and }\quad F(x,y)=0\,.
$$
For $N=72$, by proceeding similarly for the elliptic curve $36A1$ with respect its attached
 normalized newform $f\in\operatorname{New}_{36}$, we obtain the $q$-expansions of $x$ and $y$ as functions on $X_0(36)$. It is clear that the  functions for our case are $x(q^2)$ and $y(q^2)$.

\vskip 0.3 cm
\noindent{\it Step 3: Determination of a suitable generator of the extension $\Q(X_0(N))/\Q(E)$.}
Let $\cG_N$ be the multiplicative group
of the modular functions on $X_0(N)$  which are equal to $\prod_{1\leq d|N}\eta(dz)^{r_d}$  for some integers
$r_d$ and where $\eta(z)$ is the Dedekind function. The group $\cG_N$ is the multiplicative subgroup of $\Q(X_0(N))$ which consists of the normalized functions whose zeros and poles are cusps (for a detailed description of this group see 2.2 of \cite{Go91}). In our case, due to the fact that  $w$ left stable the set of cusps and  $w\in\Aut_{\Q}(X_0(N))$, the involution $w$ induces an involution $w^* $ on $\Q\otimes\cG_N$. By  Proposition 2 of \cite{Go91},  there exists
a function $u\in  \cG_N$ satisfying:
\begin{itemize}
 \item[(i)] The polar part of $u$ is a multiple of the divisor $ (\infty)$ or $(\infty)+(w(\infty))$.

 \item [(ii)] $\operatorname{div\,} w^*(u)\neq \operatorname{div\,}u$.
 \end{itemize}
Once such a function $u$ is chosen,  the divisor of $w^*(u)$ is determined and we can
find a function $v\in \cG_N$ having the same divisor as $w^*(u)$, i.e. $w^*(u)=a\,v$ for some nonzero rational number $a$. By using the $q$-expansions of $u$,  $v$, $x$ and $y$, we can determine $a$ because $u+a \,v$ must be equal to a polynomial with rational coefficients  in the  functions $x$ and $y$. In fact, if $w$ is an Atkin-Lehner involution ($N\neq 64,72$), then  Proposition 3 of \cite{Go91} allows us to determine $a$ without using $q$-expansions. In any case, with
 our choice, the rational number $a$ turns out to be always   integer. Finally, we take
then function $t:=u-a\,v$ which  satisfies $w^*(t)=-t$. Therefore, $\Q(X_0(N))=\Q(E)(t)$ and, moreover, the function $t^2$ lies in $ \Q (E)$ and has a unique pole in $\pi (\infty)$. Hence,
$$t^2=P(x,y)$$
for some polynomial $P$ with integers coefficients, which provides  an equation for $X_0(N)$ related with the chosen equation for $E$, for which its Mordell-Weil group is described in Cremona's tables.   The polynomial $P(x,y)$ is taken as a polynomial of the form $P_1(x)+P_2(x)\,y$, where
$P_1(x),P_2(x)\in\Z$. Since the degree of $P_1(x)$ agrees with $-\ord_{\infty} u$ and $\deg P_2(x)\leq \deg P_1(x)-2$, $u$ is chosen to be  $-\ord_{\infty} u$ minimal. In tables 3 and 4 of the appendix, the  corresponding functions $t$ are exhibited for the new and non-new case respectively, while   in Tables 5 and 6  the polynomials $P(x,y)$ are  presented for the new and non-new case respectively.
\newpage
\section{Appendix}

\subsection{Tables for $E$, $w$ and $h$}
$$\text{Table 1 (new case)}$$

\vskip 0.2 cm
$$
\begin{array}{c||c|rc|c}
N& w  & & X_0(N)/\langle w\rangle & h (q)\\\hline
34&w_{17}    &34A1&y^2+x y=x^3-3 x+1 &f_{34A}(q)\\ \hline
 38&w_{19}    &38B1& y^2+yx+y=x^3+x^2+1& f_{38B}(q)\\ \hline
   43&w_{43}  &43A1&y^2+y=x^3+x^2 &f_{43A}(q)\\\hline
 44&w_{11}      & 44A1&y^2=x^3+x^3+3x-1 &  f_{44A}(q)\\ \hline
 45&w_5    &45A1& y^2+yx=x^3-x^2-5&f_{45A}(q)  \\ \hline
 51&w_{51} & 51 A1& y^2+y=x^3+x^2+x-1 & f_{51A}(q)\\ \hline
53& w_{53}  &53A1& y^2+xy +y=x^3-x^2 & f_{53A}(q)\\ \hline
54&w_{27}    & 54B1& y^2+xy+y=x^3-x^2+x-1& f_{54B}(q) \\ \hline
55&w_{11}  &55A1&y^2+xy=x^3-x^2-4x+3  & f_{55A}(q)\\ \hline
56&w_7  &56A1& y^2=x^3+x+2 &f_{56A}(q) \\ \hline
 61&w_{61}  & 61A1&y^2+xy=x^3-2x+1 & f_{61A}(q) \\ \hline
62&w_{31}  &62A1& y^2+xy+y=x^3-x^2-x+1 &  f_{62A}(q)\\ \hline
64&(S_2 \cdot w_{64})^2 &64A1& y^2=x^3-4 & f_{64A}(q)\\ \hline
65& w_{65}  & 65A1&y^2+xy=x^3-x & f_{65A}(q)\\ \hline
69& w_{23}  & 69A1&y^2+xy +y=x^3-x-1  &  f_{69A}(q)\\ \hline
79& w_{79}  &79A1& y^2+xy+y=x^3+x^2-2x & f_{79A}(q)\\ \hline
83&w_{83}  &83A1&y^2+xy+y=x^3+x^2+x  &  f_{83A}(q)\\ \hline
89&w_{89}  &89A1& y^2+xy+y=x^3-x & f_{89A}(q)\\ \hline
92&w_{23}  & 92A1&y^2=x^3+x^2+2x+1  & f_{92A}(q)\\
\hline
94& w_{47}   & 94A1&y^2+xy+y=x^3-x^2-1 & f_{94A}(q) \\ \hline
101&w_{101}  &101A1&  y^2+y=x^3x^2-x-1& f_{101A}(q)\\ \hline
131& w_{131}  & 131A1& y^2+y=x^3-x^2+x &  f_{131A}(q)\\ \hline
\end{array}
$$

\vskip 0.2cm
$$\text{Table 2 (non-new case)}$$
\vskip 0.2 cm

$$
\begin{array}{c||c| c|c|c}
N& w  & & X_0(N)/\langle w\rangle& h(q) \\ \hline
 42&  w_{14}  &21A4  &  y^2+xy=x^3+x & f_{21A}(q)+2f_{21A}(q^2)\\
\hline
 60&w_{15}  &20A2&    y^2=x^3+x^2-x& f_{20A}(q)+3f_{20A}(q^3) \\ \hline
 63&w_{63}  & 21A4&  y^2+x y=x^3+x& f_{21A}(q)-3f_{21A}(q^3)\\ \hline
72&S_2  & 36A1&  y^2=x^3+1 & f_{36A}(q^2) \\ \hline
75&w_{75}  &15A8 & y^2+x\,y+y=x^3+x^2  &  f_{15A}(q)-5f_{15A}(q^5)\\ \hline
81 & w_{81}  & 27A3 &y^2+y=x^3 &f_{27A}(q)-3f_{27A}(q^3)\\ \hline
95& w_{95} &19A3 &y^2+y=x^3+x^2+x  & f_{19A}(q)-5f_{19A}(q^5) \\ \hline
119& w_{119}  &  17A4 & y^2+xy +y=x^3-x^2-x & f_{17A}(q)-7f_{17A}(q^7)\\ \hline
\end{array}
$$

\newpage
\subsection{Tables for $t$}
$$\text{Table 3 (new case)}$$

\vskip 0.2 cm
$$
\begin{array}{c||c}
N&  t \\\hline
34&\frac{\eta(2z)^4\eta(17z)^2}{\eta(z)^2\eta(34z)^4}-17\frac{\eta(z)^2\eta(34z)^4} {\eta(2z)^4\eta(17z)^2}\\ \hline
 38& \frac{\eta(2 z)^8\eta (19 z)^4}{\eta (z)^4\eta(38 z)^8}-19^2\frac{\eta (z)^4\eta(38 z)^8}{\eta(2 z)^8\eta (19 z)^4}\\ \hline
   43&\frac{\eta(z)^4}{\eta( 43 z)^4}-43^2  \frac{\eta( 43 z)^4}{\eta(z)^4} \\\hline
 44 &  \frac{\eta(4 z)^4\eta (22 z)^2}{\eta (2z)^2\eta(44z )^4}+11 \frac{\eta (2z)^2\eta(44z )^4}{\eta(4 z)^4\eta (22 z)^2}\\ \hline
 45&\frac{\eta(9z)^3\eta(15z)}{\eta(3z)\eta(45z)^3}+5 \frac{\eta(3z)\eta(45z)^3} {\eta(9z)^3\eta(15z)}\\ \hline
51& \frac{\eta(3z)^9\eta(17z)^3}{\eta(z)^3\eta(51 z)^9} +17^3\frac{\eta(z)^3\eta(51 z)^9}{\eta(3z)^9\eta(17z)^3}\\ \hline
53& \frac{\eta(z)^6}{\eta(53z)^6}-53^3 \frac{\eta(53 z)^6}{\eta(z)^6}\\ \hline
54& \frac{ \eta(18z)\eta(27 z)^3}{\eta(9z)\eta(54z)^3}+\frac{ \eta(z)^3\eta(6z)}{\eta(2z)^3\eta(3z)} \\ \hline
55& \frac{\eta(5z)^5\eta(11z)}{\eta( z)\eta(55 z)^5}-11^2\frac{\eta( z)\eta(55 z)^5}{\eta(5z)^5\eta(11z)}\\ \hline
56 &\frac{\eta(8z)^4\eta(28 z)^2}{\eta(4z)^2\eta(56 z)^4}+7 \frac{\eta(4z)^2\eta(56 z)^4}{\eta(8z)^4\eta(28 z)^2} \\ \hline
 61& \frac{\eta (z)^2}{\eta(61z)^2}-61 \frac{\eta(61z)^2}{\eta (z)^2} \\ \hline
62 &  \frac{\eta(2z)^8\eta(31 z)^4}{\eta(z)^4\eta(62 z)^8}-31^2\, \frac{\eta(z)^4\eta(62 z)^8}{\eta(2z)^8\eta(31 z)^4}\\ \hline
64 & \frac{\eta(32z)^6}{\eta(16z)^2\eta(64 z)^4}-4 \frac{\eta(16z)^2\eta(64 z)^4}{\eta(32z)^6}\\ \hline
65 & \frac{\eta (5z)^5\eta(13 z)}{\eta (z)\eta(65z)^5} -13^2 \frac{ \eta (5z) \eta(13 z)^5}{\eta(z)^5\eta (65 z)}\\ \hline
69  &  \frac{\eta(3z)^9\eta (23 z)^3}{\eta(z)^3\eta(69 z )^9}+23^3 \frac{\eta(z)^3\eta(69z )^9}{\eta(3z)^9\eta (23 z)^3}\\ \hline
79 & \frac{\eta (z)^4}{\eta(79z)^4}- 79^2\frac{\eta(79 z)^4}{\eta( z)^4}\\ \hline
83 &  \frac{\eta{z}^12}{\eta(83 z)^12}-83^6\, \frac{\eta(83 z)^12}{\eta{z}^12}\\ \hline
89 &  \frac{\eta(z)^6}{\eta(89z)^6}-89^3 \frac{\eta(89z)^6}{\eta(z)^6}\\ \hline
92  & \frac{\eta (4z)^4\eta(46z)^2}{\eta(2z)^2\eta(92z)^4}+23 \frac{\eta(2z)^2\eta(92z)^4}{\eta (4z)^4\eta(46z)^2}\\
\hline
94 & \frac{\eta(2z)^8\eta(47z)^4}{\eta(z)^3\eta(94 z)^8}-47^2\, \frac{\eta(z)^3\eta(94 z)^8}{\eta(2z)^8\eta(47z)^4} \\ \hline
101&  \frac{\eta(z)^6}{\eta(101 z)^6}-101^3\, \frac{\eta(101 z)^6}{\eta(z)^6}\\ \hline
131&  \frac{\eta(z)^{12}}{\eta(131z)^{12}}-131^6\frac{\eta(131z)^{12}}{\eta(z)^{12}}\\ \hline
\end{array}
$$

\skip 0.4 cm
$$\text{Table 4 (non-new case)}$$

\skip 0.2 cm
$$
\begin{array}{c||c}
N&  t \\\hline
 42 & \frac{\eta (z)^9 \eta(2z)^3\eta(6z)^5 \eta(14z)^3\eta(21z)^7}{\eta(3z)^{13}\eta(7z)^3\eta (42 z)^{11}} -7^2\frac{\eta (z)^3 \eta(6z)^7\eta(7z)^3 \eta(14z)^9\eta(21z)^5}{\eta(2z)^{3}\eta(3z)^{11}\eta (42 z)^{13}}   \\ \hline
 60&  \frac{\eta(2z) \eta(12 z)^6 \eta(20 z)^2\eta(30 z)^3}{\eta( 4 z)^2\eta(6 z)^3\eta(10z)\eta(60z)^6}-5  \frac{\eta(2 z)^3\eta(12 z)^2\eta (20 z)^6 \eta(30 z)}{\eta(4z)^6 \eta(6 z)\eta(10 z)^3\eta(60 z)^2}\\ \hline
 63& \frac{\eta(9z)^3\eta(21 z)}{\eta(3z)\eta(63z)^3}-7 \frac{ \eta(3z) \eta(7z)^3}{\eta z)^3\eta(21z)}\\ \hline
72 &  \frac{ \eta(z)^6 \eta(6 z)\eta( 24 z)^2\eta(36 z)^3}{\eta (2 z)^3\eta(3z )^2\eta (12 z)\eta(72 z)^6}-
\frac{ \eta(2z)^{15} \eta(3 z)^2 \eta( 12 z) \eta(24 z)^2\eta (36 z)^3}{\eta ( z)^6\eta(4z )^6\eta (6 z)^5 \eta(72 z)^6}\\ \hline
75 & \frac{\eta (3z)^3\eta (25 z)}{\eta(z)\eta(75z)^3}-5^2\frac{\eta 3z) \eta(25 z)^3}{\eta z)^3\eta(75 z)} \\ \hline
81  &\frac{\eta (z)^3 \eta(27 z)}{\eta (3 z) \eta(81 z)^3}-3^5\frac{\eta (3 z) \eta(81 z)^3}{\eta (z)^3 \eta(27 z)}\\ \hline
95 &  \frac{\eta(5z)^5\eta (19z)}{\eta(z)\eta(95 z)^5}-19^2\frac{\eta(5 z) \eta (19 z)^5}{\eta (z)^5\eta (95 z)}\\ \hline
119 &   \frac{\eta( 7z)^7\eta(17 z)}{\eta z)\eta(119 z)^7}-17^3 \frac{\eta (7z)\eta(17 z)^7}{\eta z)^7\eta(119 z)}\\ \hline
\end{array}
$$
\newpage
\subsection{Tables for $P(x,y)$}

$$ \text{Table 5 (new case)}$$
$$
\begin{array}{c|l}
N&  P(x,y)\\\hline
34 & -48 - 32x + 20x^2 + 24x^3 + x^4 +8 ( 2+2x+x^2)y \\
\hline
 38  & -960 + 3168x + 13160x^2 + 21724x^3 + 25833x^4 + 21810x^5 + 10071x^6 + 2065x^7 + 136x^8 + x^9 \\
  & +x(56 + 27x + x^2)(44 + 88x + 137x^2 + 102x^3 + 17x^4)y \\ \hline
  43 & -7200 - 1680x + 9400x^2 - 2332x^3 - 4868x^4 + 1708x^5 + 194x^6 + x^7\\
  &  -(-72 + 36x + x^2)(-44 - 84x + 45x^2 + 22x^3)y\\\hline
   44 &(2 + x)^2(7 + 3x + 5x^2 + x^3)\\ \hline
   45 &  x^2(-3 + 6x + x^2 + 4y)  \\ \hline
   51 &-8904 + 89496x + 720815x^2 + 2136731x^3 + 3806784x^4 + 4786996x^5 + 4564407x^6 + 3440158x^7  \\
    & + 2089704x^8 + 1029855x^9 + 409276x^{10 }+ 129052x^{11} +    31311x^{12} + 5557x^{13} + 658x^{14} + 43x^{15} + x^{16}  \\
     &   +(-5 + 3x + 2x^2 + x^3)(36 + 37x + 20x^2 + 3x^3)  \\
     &  (304 + 1445x + 2641x^2 + 2567x^3 + 1636x^4 + 706x^5 +   201x^6 + 34x^7 + 2x^8) y   \\ \hline
     53 &  -247408 + 665520x - 1831348x^2 + 4346036x^3 - 7515167x^4 + 7342874x^5 - 4503204x^6  \\
      & + 2095505x^7 - 818846x^8 + 230692x^9 - 33955x^{10} + 739x^{11} +    237x^{12} + x^{13} -(148 -   1108x \\
      & + 151x^2 + 1363x^3 - 621x^4 - 12x^5 + 25x^6) (1328 - 1308x + 805x^2 - 328x^3 + 47x^4 + x^5) y\\
      \hline
       54 & 3 - 3x + 3x^2 + x^3+ 3(1 + x) y \\ \hline
       55 & (2 + x)(-138 - 271x - 58x^2 + 1411x^3 + 168x^4 - 1461x^5 -   281x^6 + 349x^7 + 68x^8 + x^9) \\
       &  + (2 + x)(-1 + 2x)(4 + 3x)(8 - 26x - 78x^2 + 11x^3 + 28x^4 +   2x^5)y  \\\hline
       56 &(7 + x^2)(2 - x + x^2)(2 + x + x^2)\\ \hline
       61 &-122 + 176x - 27x^2 - 65x^3 + 18x^4 + x^5-(1 + x)(-22 - 15x + 9x^2)y\\ \hline
        62 & -3840 - 448x + 12724x^2 + 42628x^3 + 62861x^4 + 5174x^5 +  109639x^6 + 289900x^7 +73179x^8 \\
        &  - 61722x^9 + 143262x^{10} +  178641x^{11} + 61858x^{12 }+     7490x^{13} + 253x^{14} + x^{15}+xy (-4 + 192x\\
        &  +    446x^2 - 108x^3 - 268x^4 + 443x^5 + 284x^6 + 23x^7)(32 - 56x + 178x^2 + 248x^3 + 55x^4 + x^5)\\
        \hline
         64 & (-2 + x)(2 + x)(4 + x^2)\\\hline
         65 & 1 - 35x - 85x^2 - 15x^3 - 35x^4 - 50x^5 - 403x^6 + 10x^7 +  35x^8 + 65x^9 - 90x^{10 }+ 25x^{11} + x^{12} \\
         & -5(-1 + x)(1 + x)(1 + x^2)(-1 - 2x + x^2)(2 - 5x + 2x^2 +    5x^3 + 2x^4)y \\\hline
          69 &40128 - 2804032x - 24658412x^2 - 82258148x^3 - 78001407x^4 + 286063638x^5 + 1082537261x^6   \\
          & + 1420597832x^7+ 9621058x^8 -2694511846x^9 -    4047900698x^{10} - 2330523372x^{11}  + 840632694x^{12}\\
          & +  2638911745x^{13 }+ 2331856822x^{14 }+ 1199425309x^{15 }+393442428x^{16}  + 82393205x^{17}+ 10602593x^{18}\\
          & +     779531x^{19} + 28810x^{20} + 412x^{21 }+ x^{22}+y (-432 - 2920x - 8036x^2 - 4860x^3 + 11574x^4+ 21734x^5\\
          &  +   14665x^6 + 4288x^7 + 477x^8 + 14x^9)  (1276 + 3128x - 4870x^2 - 28854x^3 - 3924x^4 - 657x^5 +   3703x^6 \\
                     & + 40883x^7 + 1804x^8 + 3299x^9 + 202x^{10}+ 2x^{11}) \\ \hline
           79& -24843 - 7420x + 112556x^2 + 76149x^3 - 214447x^4 - 113728x^5 +  157812x^6 + 73431x^7 \\
           & - 49467x^8 - 22769x^9 + 5008x^{10 }+  2736x^{11} + 181x^{12} + x^{13} \\
           & -(1 + x)(83 - 145x - 10x^2 + 39x^3 + x^4)(105 + 599x -    75x^2 - 604x^3 - 31x^4 + 141x^5 + 21x^6)y  \\\hline
            83 & -846820980000 + 1701842643824x - 4190038951864x^2 - 15407944317740x^3 + 52631374705524x^4  \\
             &+ 195902048285636x^5-69755046878975x^6 -    1014877154551415x^7 - 1063602170418749x^8 \\
             &+ 1855157981145929x^9   + 5075380899888979x^{10} + 1636529117010692x^{11} - 8302874421713802x^{12}         \\
             & -    11678206852543817x^{13} + 1005402492172935x^{14} + 17071231491541350x^{15} + 14609656638884595 x^{16}   \\
             &- 5861471722333698 x^{17} - 19444452135637043x^{18 }- 10558933244522770x^{19 }+ 7380101893789387x^{20} \\
             & + 13072119010688686x^{21} + 4488232204563914x^{22} -4533961869101651x^{23} - 5251219501592566x^{24}\\
               &  -    1152837411146407x^{25} + 1460566880152011x^{26 }+1235263008419117x^{27 }+ 207233373480590x^{28}\\
               & - 227667918937852x^{29} - 159972030244333x^{30} - 31212179742782x^{31} + 11503837595608x^{32}   \\
               & + 9040609302177x^{33} + 2671633081498x^{34} + 434275870731x^{35}+ 39994913022x^{36} + 2002817221x^{37}\\
               & + 50297678x^{38} +    551168x^{39} + 1991x^{40} + x^{41 } -y(452104 - 2937060x - 9114834x^2+ 6837837x^3 \\
                 &  +     31676870x^4 + 15747540x^5 - 39815725x^6- 52447587x^7+3221373x^8 + 46993123x^9 + 26434445x^{10}\\
                 &   - 10569436x^{11} -     16815463x^{12} - 4439206x^{13} + 2240588x^{14} + 1704048x^{15} +     372443x^{16 }+ 26456x^{17}  \\ & +  483x^{18} + x^{19})    (1809956 - 279348x - 16404210x^2 - 14857887x^3 +   44872703x^4 + 91386412x^5  \\
                    & + 7209319x^6 - 135817396x^7 -134605986x^8 + 24499468x^9 +       132852372x^{10}+ 79321327x^{11}  \\
                      &  - 23111323x^{12} -   50396693x^{13} - 19474566x^{14} + 4349474x^{15}  + 6065443x^{16} +   1970373x^{17}+236111x^{18 }  \\
                       &+ 8734x^{19} +       65x^{20})\\ \hline
 \end{array}
 $$

  \newpage
$$
\begin{array}{c|l}\hline
89 &  -17600 + 410400x - 27548480x^2 + 10400948x^3 + 146498188x^4 -  32027037x^5 - 360910680x^6  \\
&  - 17199072x^7 + 501894798x^8 +  161240831x^9 - 391130731x^{10} -     233566274x^{11 }+ 145795788x^{12} \\
& + 151162884x^{13} - 4448240x^{14} - 44442246x^{15} - 12762946x^{16} + 3535663x^{17}+  2568284x^{18}\\
&  + 483388x^{19} + 30642x^{20} +     515x^{21} +  x^{22} - y(-1 + x)(1 + x)(2700 + 3240x - 3609x^2- 4873x^3\\
&   +    323x^4 + 1473x^5 + 465x^6 + 17x^7)(648 - 1340x - 10451x^2 + 10324x^3 + 23063x^4 - 2172x^5 \\
& -    17100x^6- 5730x^7  + 3597x^8 + 2256x^9 + 243x^{10} + 2x^{11}) \\ \hline
 92 & (3 + x)^2(8 + 4x + 4x^2 + x^3)(8 + 20x + 28x^2 + 25x^3 +
   14x^4 + 5x^5 + x^6)\\ \hline
    94 & -21715 - 47508 x - 103195 x^2 + 219398 x^3 + 1663909 x^4 +  5469799 x^5 + 10685097 x^6  \\
    &   + 13118353 x^7 + 4598983 x^8 -  18554364 x^9 - 49262084 x^{10} -     67097732 x^{11} - 54688267 x^{12}    \\
    &- 12748826 x^{13 }   + 34980862 x^{14 }+  61712870 x^{15} + 58008344 x^{16} + 36328896 x^{17} + 15246919 x^{18} \\
    &  +  3869878 x^{19 }+ 504954 x^{20} +     27667 x^{21} + 465 x^{22} +  x^{23 } +y(-120 - 313 x - 845 x^2 - 566 x^3    \\
    & + 598 x^4 + 2600 x^5 +    3476 x^6 + 2830 x^7 + 1066 x^8 + 105 x^9 + x^{10}) (42 - 275 x - 1064 x^2 - 2649 x^3 \\
    &   - 3694 x^4 - 2814 x^5 +    602 x^6+ 4402 x^7 + 5384 x^8 + 3641 x^9 + 814 x^{10} + 31 x^{11}) \\
   \hline
    101 &  3973376 + 24345712 x - 57185964 x^2 - 178101464 x^3 +  161167773 x^4 + 462091312 x^5 \\
    & - 328497531 x^6- 579757947 x^7 +  504756351 x^8 + 320809112 x^9 -     453760989 x^{10} + 7488095 x^{11} \\
    & + 195974239 x^{12} - 88118541 x^{13} -  19312382 x^{14 }+ 31143884 x^{15} - 10685229 x^{16} - 391854 x^{17 } \\
    & +  1833695 x^{18}- 900817 x^{19} +     273859 x^{20} - 59884 x^{21 }+ 9221 x^{22} - 857 x^{23} + 29 x^{24 }+  x^{25}\\
        & -y(-2 + x)(-2492 - 5984 x +    5085 x^2 + 12629 x^3 - 5462 x^4 - 7457 x^5 + 4351 x^6+ 561 x^7  \\ & -    1100 x^8+ 476 x^9-        115 x^{10 } + 12 x^{11})(-1164 + 4596 x + 5916 x^2 - 11845 x^3 +   1691 x^4 + 5276 x^5 \\
    & - 3587 x^6 +        849 x^7- 17 x^8 - 12 x^9 - 4 x^{10} + x^{11})\\ \hline
    131 &-11491793287200 - 77827916513200x + 271274799348200x^2 +  2891920043667900x^3   \\
     & - 7413871118879700x^4-208181329634260700x^5 + 1446446595356863725x^6 -     676025577299157695x^7 \\
     &  - 34278119003227037765x^8 +  204438101122657660450x^9 - 568517654861536891120x^{10} \\
     &  +  481549779729164516941x^{11}+2635245766368675205851x^{12} - 13175502704057926427955x^{13} \\
      & +  31767041079352267060729x^{14} - 43173166155030136041073x^{15} +
 8969030371497162128939x^{16} \\
     &  +     112429758234629584481836x^{17} - 310370654615807411016501x^{18} +470937641418992597885012x^{19} \\
     & - 405839688122220307045126x^{20 }-  8033748795341970931271x^{21} +     668619555576572629194039x^{22 } \\
      & - 1234109512950759734621090x^{23} +  1327319462648814369990544x^{24}   \\
     & - 822566308055883209677805x^{25}-33736964013542515967889x^{26} + 788027884980919155268540x^{27}   \\
     & -  1093171730235629694196110x^{28}+ 911004881546511370059093x^{29} -     467554555186095909777550x^{30}   \\
       & + 47296265160342002570173x^{31}+  185250143048432609363127x^{32} - 225960509043141561845001x^{33} \\
       & + 160852633886373809385204x^{34 }-     77496232791333991087714x^{35} + 20802958115704540858490x^{36}  \\
       & + 3880011388402181168118x^{37}  - 8735192594907803099043x^{38} +  6210755007512712251897x^{39} \\
       & -    3003413163319348531225x^{40}  + 1070020657667042767223 x^{41} -   259403449388540446685 x^{42}     \\
       &  +20958569399527863492 x^{44} -     14326879082941335638 x^{45} + 5912238306538064916 x^{46 } \\
       & + 18477953921830947985 x^{43} -  1877971630289666366 x^{47 }+ 488283858739094050 x^{48}   \\
       &-  105722203629934910 x^{49}+ 18923329741114761 x^{50} -     2689272320340289 x^{51} + 263300065498249 x^{52}
       \\
   &  - 4260505636661 x^{53} - 5157557959005 x^{54} + 1416554294641 x^{55} -  243661620624 x^{56} + 31240505016 x^{57} \\
       & - 2997027299 x^{58} + 189957325 x^{59} - 2699879 x^{60} - 1016956 x^{61} +  130698 x^{62} - 7747 x^{63} + 187 x^{64} + x^{65} \\
       &  -y(-1 + x)(-1651336 - 45600164 x - 27727502 x^2 +    1818599597 x^3 - 7927670906 x^4 \\
       & + 13960778588 x^5-421393207 x^6- 51398069731 x^7 + 119271135931 x^8 -        134713281163 x^9  \\
       & + 52661398574 x^{10} +80481104397 x^{11} -    163994123959 x^{12} + 148050275877 x^{13} - 72614888314 x^{14}  \\
         &+    6183551272 x^{15} + 20717235269 x^{16}  -  19147906676 x^{17} + 9847397087 x^{18} - 3367121265 x^{19} \\
          & +     730146720 x^{20 }- 51809318 x^{21} - 30959822 x^{22} + 15666228 x^{23} - 4292594 x^{24} + 835815 x^{25} -        122075 x^{26}
           \\
            & + 13084 x^{27} - 903 x^{28} + 23 x^{29} +    x^{30}) (-4255900 + 71222100 x + 1461150 x^2  - 2150436025 x^3 \\
              & +    7702741575 x^4- 2543947440 x^5 -421393207 x^6  - 51398069731 x^7+ 119271135931 x^8\\
              & -        134713281163 x^9  + 52661398574 x^{10} + 80481104397 x^{11 }-163994123959 x^{12} + 148050275877 x^{13} \\
              & - 72614888314 x^{14} +    6183551272 x^{15} + 20717235269 x^{16}-    19147906676 x^{17} + 9847397087 x^{18}  \\
                 &- 3367121265 x^{19}+     730146720 x^{20 }- 51809318 x^{21} - 30959822 x^{22}+ 15666228 x^{23} - 4292594 x^{24}  \\
                  & + 835815 x^{25} -        122075 x^{26}+ 13084 x^{27} - 903 x^{28} + 23 x^{29} +    x^{30})(-4255900 + 71222100 x + 1461150 x^2  \\
                  & - 2150436025 x^3+    7702741575 x^4 - 2543947440 x^5 -
        55631723695 x^6 + 194169167040 x^7   \\
  & - 323648869980 x^8  +  226360273540 x^9 + 229278613251 x^{10} - 835109705023 x^{11} +  1135514398082 x^{12 }\\
   & - 867152409818 x^{13 }+        237463102523 x^{14} + 300544359512 x^{15 }- 479908400891 x^{16} +368024466511 x^{17 }\\
  & - 177471360585 x^{18 } + 47572071000 x^{19} +    3135089681 x^{20} - 10639183810 x^{21} +    6239233500 x^{22}  \\
    & - 2354646874 x^{23}+ 659880730 x^{24 }-    144046879 x^{25} + 25143814 x^{26} - 3620625 x^{27}+ 455809 x^{28}  \\ &- 54178 x^{29} + 6052 x^{30}- 527 x^{31} +        24 x^{32}) \\\hline
\end{array}
$$

\newpage

\vskip 0.3 cm

 $$\text{Table 6 (non-new case)}$$
 \vskip 0.2 cm
 $$
\begin{array}{c|l}
N&  P(x,y)\\\hline
 42& 9(4 + x + 4x^2)\left (64 + 1017x + 96x^2 + 1178x^3 - 1352x^4 + 2883x^5 - 1336x^6 +
 730x^7 + 1800x^8 \right.\\
  & \left. + 1417x^9 + 64x^{10}+ 72 y(-1 + x)(1 + x)(50 + 213x - 6x^2 + 215x^3 -
   6x^4 + 213x^5 + 50x^6) \right)\\\hline
  60 &(-1 + x + x^2)(-1 + 4x + x^2)(1 - x + 2x^2 + x^3 + x^4)\\ \hline
  63 &1 + 9x - 24x^2 - 8x^3 - 9x^4 + 3x^5 +  x^6 -3(-2 + x)(-1 + x)(1 + x)(-1 + 2x)y\\ \hline
  72 & 4(7 + 144x + 72x^2 + 72x^3 + 144x^4 + 72x^5)+144(1 + x)(-2 +
   6x + x^3) y \\
 \hline
   75 &           -99 - 390x - 569x^2 - 372x^3 - 67x^4 + 56x^5 + 40x^6 + 11x^7                  \\
     &    -3x(2 + x)(-2 - 4x + 3x^2 + 5x^3 + 2x^4) y\\ \hline
     81 &   -968 + 132x + 837x^2 - 3213x^3 + 4107x^4 - 2223x^5 + 510x^6 -
 21x^7 - 9x^8 +
 x^9         \\
      &   -3(-1 + x)(-29 + 27x - 9x^2 + x^3)(11 + 24x - 15x^2 +
   2x^3)       y     \\ \hline
   95 & x(1 - 1208x - 10934x^2 - 44162x^3 - 109477x^4 - 180353x^5 -
    196536x^6 - 134741x^7 - 40197x^8 + \\
   & 17286x^9 +    20552x^{10} +
    4156x^{11} - 2604x^{12 }-
        1110x^{13} + 370x^{14} + 18x^{15} - 13 x^{16} + x^{17}) \\
          & + (1 +
    x)(-1 - 7x - 1218x^2 - 8470x^3 - 25928x^4 - 48038x^5 -
    51018x^6 - 26540x^7 + \\
                    & 2089x^8 +
        10997x^9 + 2802x^{10} - 1822 x^{11} - 545 x^{12}+ 323x^{13 }-
    46 x^{14} + 2x^{15}) y \\\hline
    119 & -19456 + 128127x + 31684x^2 - 1935597x^3 + 2686286x^4 +
 12402399x^5 - 29855351x^6 - 41030815x^7 \\
      & + 159290916x^8 +
 60745925x^9 - 524982545x^{10}+ 38052811x^{11} + 1173318320x^{12} - 347340496x^{13} \\
      &-1864043953x^{14} + 721098268x^{15} + 2154666360 x^{16} -
 836975880x^{17} - 1819472378x^{18}  \\
      &+ 598982372 x^{19}+
    1104865348 x^{20} - 257570193x^{21} - 464230322x^{22} +
 55914351x^{23} + 126047469x^{24} \\
    & - 504853x^{25} - 19453402x^{26} -
 2191709x^{27} + 1202359x^{28} + 297037 x^{29} +
    20766 x^{30} + 407 x^{31} \\
    &-(1 + x)(19457 - 147570x + 135420x^2 + 1672260x^3 -
   4351085x^4 - 6352626x^5 + 33446252x^6\\
    & - 1612436x^7 -
   131903527x^8+93891838x^9 + 312797860x^{10} - 358809060x^{11} -
   481203790x^{12}  \\
    &+ 724102984x^{13}+ 508610138x^{14} -
   917130440x^{15} - 392684774x^{16} + 764808468x^{17} +
       236037174x^{18}\\
      & - 417563192x^{19} - 113096097x^{20} +
   141052022x^{21} + 40143942x^{22} - 25553396x^{23} \\
        &- 8853365x^{24} +
   1504686 x^{25 }+ 853274x^{26} + 100016x^{27} +
       3681x^{28} + 30x^{29})y\\
     \hline
  \end{array}
$$

\bibliographystyle{plain}
\bibliography{bielliptic}

\end{document}